\documentclass[12pt,a4paper]{amsart}
\usepackage{amssymb,amsmath,amsthm}
\usepackage{graphicx}
\usepackage{hyperref}
 \usepackage[color=green!40]{todonotes}
\usepackage{mathtools}
\usepackage{enumerate}

\newcommand\cF{{\mathcal F}}

\newcommand\cG{{\mathcal G}}

\newcommand\cT{{\mathcal T}}

\makeatletter
\newtheorem*{rep@theorem}{\rep@title}
\newcommand{\newreptheorem}[2]{%
\newenvironment{rep#1}[1]{%
 \def\rep@title{#2 \ref{##1}}%
 \begin{rep@theorem}}%
 {\end{rep@theorem}}}
\makeatother

\theoremstyle{plain}
\newtheorem{theorem}{Theorem}[section]
\newreptheorem{theorem}{Theorem}
\newtheorem{lemma}[theorem]{Lemma}
\newtheorem{corollary}[theorem]{Corollary}

\theoremstyle{definition}

\newtheorem{defn}[theorem]{Definition}

\newcommand\cref[1]{Corollary~\ref{cor:#1}}

\title{Size, diversity, minimum degree, sturdiness, dömdödöm}

\date{}
\author{Bal\'azs Patk\'os} 
\address{HUN-REN Alfr\'ed R\'enyi Institute of Mathematics} 
\email{patkos@renyi.hu}
\begin{document}

\begin{abstract}
For a family $\cF$ of sets and a disjoint pair $A,B$ we let $\cF(A,\overline{B})=\{F\in \cF: A\subseteq F, ~B\cap F=\emptyset\}$. The \textbf{$(p,q)$-dömdödöm} of a family $\cF\subseteq 2^{[n]}$ is $\beta_{p,q}(\cF)=\min\{|\cF(A,\overline{B})|:|A|=p,|B|=q, A\cap B=\emptyset, A,B\subseteq [n]\} $. This definition encompasses size, diversity, minimum degree, and sturdiness as special cases. We investigate the maximum possible  value $\beta_{p,q}(n,k)$ of $\beta_{p,q}(\cF)$ over all $k$-uniform intersecting families $\cF\subset 2^{[n]}$. We determine the order of magnitude of $\beta_{p,q}(n,k)$ for all fixed $p,q,k$. We relate the asymptotics of $\beta_{p,q}(n,k)$ to the constant value of $\beta_{0,q}(n,q+1)$ and establish $\beta_{p,1}(n,k)=\binom{n-3-p}{k-2-p}$ and $\beta_{p,2}(n,k)=2\binom{n-5}{k-3-p}-\binom{n-7}{k-5-p}$ if $n$ is large enough.    
\end{abstract}
\maketitle

\section{Introduction}

We use standard notation: $[n]$ denotes the set of the first $n$ positive integers. For a set $X$, we use $2^X$ and $\binom{X}{k}$ to denote its power set and the family of all its $k$-subsets, respectively.

In this short note, we address extremal problems on intersecting families of sets. A family $\cF$ of sets is \textbf{intersecting} if $F\cap F'\neq \emptyset$ holds for all $F,F'\in \cF$. Our main definition is as follows.

\begin{defn}
    For a family $\cF$ of sets and a disjoint pair $A,B$ we let $$\cF(A,\overline{B})=\{F\in \cF: A\subseteq F, ~B\cap F=\emptyset\}$$
    and if $A=\{i\}$ or $B=\{j\}$, we write $\cF(i,\overline{B})$ and $\cF(A,\overline{j})$. The \textbf{$(p,q)$-dömdödöm}\footnote{Dömdödöm is a Hungarian word that does not mean anything, but can mean anything based on the speaker's intent, and yet is unknown even to most Hungarians. Dömdödöm is one of the main characters of Ervin Lázár's Square-shaped Round Forest \cite{L}. This somewhat depressed creature refuses to say anything other than dömdödöm, and it is only his friends who can tell whether he meant "I don't like this" or "I love you all".} of a family $\cF\subseteq 2^{[n]}$ is $$\beta_{p,q}(\cF)=\min\{|\cF(A,\overline{B})|:|A|=p,|B|=q, A\cap B=\emptyset, A,B\subseteq [n]\}.$$
\end{defn} 
Observe that $\beta_{0,0}(\cF)=|\cF|$, $\beta_{1,0}(\cF)$ is the minimum degree $\delta(\cF)$, $\beta_{0,1}(\cF)$ is $|\cF|-\Delta(\cF)=\gamma(\cF)$ the diversity (formerly known as the unbalance \cite{LP}) of $\cF$, and $\beta_{1,1}(\cF)$ is the sturdiness \cite{FW} of $\cF$.

We will be interested in the maximum possible value of $\beta_{p,q}(\cF)$ over all $k$-uniform intersecting families $\cF\subset 2^{[n]}$. Formally, we define
\[
\beta_{p,q}(n,k):=\max\left\{\beta_{p,q}(\cF): \cF\subseteq \binom{[n]}{k}, \cF ~\text{is intersecting}\right\}.
\]

The statement $\beta_{0,0}(n,k)=\binom{n-1}{k-1}$ for $n\ge 2k$ is the celebrated theorem of Erd\H os, Ko, and Rado \cite{EKR}.

$\beta_{0,1}(n,k)=\binom{n-3}{k-2}$ for large enough $n$ was first established in \cite{LP}, and then the threshold for $n$ was improved several times \cite{F,Kup} and the current best bound $n\ge 36k$ is due to Frankl and Wang \cite{FWdiv}.

$\beta_{1,0}(n,k)=\binom{n-2}{k-2}$ states that maximum possible minimum degree in a $k$-uniform intersecting family is achieved by the star if $n>2k$ and was proved by Huang and Zhao \cite{HZ}. Later Kupavskii determined \cite{Kupdeg} the maximum of the minimum $t$-degree of an intersecting family, i.e. established $\beta_{q,0}(n,k)=\binom{n-1-q}{k-1-q}$ if $n\ge 2k+ \frac{3q}{1-\frac{t}{k}}$.

Finally, the maximum sturdiness of intersecting families was determined by Frankl and Wang \cite{FW}: $\beta_{1,1}(n,k)=\binom{n-4}{k-3}$ if $n\ge 36(k+6)$.

The \textit{covering number} $\tau(\cF)$ is the minimum size of a set $X$ such that $X\cap F$ is non-empty for every $F\in \cF$. Sets with this property are called \textit{transversals} or \textit{covers} of $\cF$. A \textit{minimal cover} is a cover $X$ such that any proper subset $Y\subsetneq X$ is not a cover. The family of minimal covers of $\cF$ is denoted by $\cT_\cF$. The relationship of the covering number and the size of an intersecting family has been well studied \cite{Fcov,FKlarge,FOT}. Our first result relates the $(p,q)$-dömdödöm of an intersecting family to its covering number.

\begin{theorem}\label{tau}
    Let $k,p,q$ be fixed and $\cF=\cF_n\subseteq \binom{[n]}{k}$ an intersecting family. If $\tau(\cF)\le q$, then $\beta_{p,q}(\cF)=0$. If $\tau(\cF)>q$, then $\beta_{p,q}(\cF)=O(n^{k-p-\tau(\cF)})$.
\end{theorem}

As a corollary, we almost immediately obtain the order of magnitude of $\beta_{p,q}(n,k)$ for all fixed $k,p,q$ and $n$ tending to infinity.
\begin{corollary}\label{order}
    For fixed $p,q,k$ the quantity $\beta_{p,q}(n,k)$ is positive for all large enough $n$ if and only if $p+q<k$ holds. Moreover, if $p+q<k$, then $\beta_{p,q}(n,k)=\Theta(n^{k-1-p-q})$.
\end{corollary}


It is known \cite{EL,T} that the number of vertices in $k$-uniform intersecting families with covering number $k$ is bounded by a constant $f(k)$. Therefore, by Theorem \ref{tau}, we have $\beta_{0,q}(n',q+1)=\beta_{0,q}(n'',q+1)$ if both $n'$ and $n''$ are large enough. We denote this constant value by $\beta(q)$. Our next result states that $\beta(q)$ determines the asymptotics of $\beta_{p,q}(n,k)$ for all $k\ge p+q+1$.

\begin{theorem}\label{0q}
    For any $q\ge 1$ and $k\ge q+1$ fixed, we have $\beta_{p,q}(n,k)=(\beta(q)+o(1))\binom{n-2q-p-1}{k-p-q-1}$.
\end{theorem}

Finally, we determine $\beta_{p,1}(n,k)$ and $\beta_{p,2}(n,k)$ for all $k\ge 3$ and $n\ge n_0(k,p)$.

\begin{theorem}\label{02}
\
    \begin{enumerate}
        \item 
        For any $k\ge p+2$ there exists $n_0(k,p)$  such that we have $\beta_{p,1}(n,k)=\binom{n-3-p}{k-2-p}$ if $n\ge n_0(k,p)$.
        \item 
        For any  $k\ge p+3$ there exists $n_0(k,p)$ such that we have $\beta_{p,2}(n,k)=2\binom{n-5-p}{k-3-p}-\binom{n-7-p}{k-5-p}$ if $n\ge n_0(k,p)$.
    \end{enumerate} 
\end{theorem}

\section{Proofs}

Before we start proving our results, let us gather some important observations that we will use frequently in our arguments. First, by definition of $\beta_{p,q}$, if $\cG\subseteq \cG'$, then $\beta_{p,q}(\cG)\le \beta_{p,q}(\cG')$ for any $p,q\ge 0$. So in all our upper bound proofs we will assume $\cF$ is maximal intersecting. Let $\cT_\cF$ be the family of minimal covers of $\cF$. Then, as $\cF$ is maximal intersecting, we have $\cF=\cup_{T\in \cT_\cF}\nabla^n_k(T)$, where $\nabla^n_k(T)=\{G\in \binom{[n]}{k}:T\subseteq G\}$. Also, $\cT_\cF$ is intersecting as if $T,T'$ are disjoint members of $\cT_\cF$, then $\nabla^n_k(T)$ and $\nabla^n_k(T')$ contain a pair of disjoint sets contradicting the intersecting property of $\cF$. We will also need the following result.

\begin{theorem}[\cite{EL,Gy,T}]\label{cover}
    The number of minimal covers of size at most $k$ of a $k$-uniform family is at most $k^k$.
\end{theorem}

\begin{proof}[\textbf{Proof of Theorem \ref{tau}}]
    Let $\cF\subseteq \binom{[n]}{k}$ be a $k$-uniform intersecting family.  If $\tau(\cF)\le q$, then consider any $q$-set $Q$ that contains a cover $T\in \cT_\cF$. Then, as $T$ is a cover, there is no $F\in \cF$ that is disjoint with $T$ let alone with $Q$ and thus $\beta_{p,q}(\cF)=0$. Assume next that $\tau(\cT)>q$. By Theorem \ref{cover}, $|\cup_{T\in \cT_\cF}T|\le k^{k+1}$. So if $n\ge k^{k+1}+p$, then there exists a $p$-subset $P$ of $[n]$ such that $P\cap (\cup_{T\in \cT_\cF}T)=\emptyset$. Any member $F$ of $\cF$ containing $P$ should contain $T\cup P$ for some $T\in \cT_\cF$. But then as any $T\in \cT_\cF$ has size at least $\cT(\cF)$, $P \cup T$ has size at least $p+\tau(\cF)$, and so $\cF(P,\overline{Q})\le k^{k+1}\binom{n-p-\tau(\cF)}{k-p-\tau(\cF)}=O(n^{k-p-\tau(\cF)})$ as claimed.
\end{proof}

\begin{proof}[\textbf{Proof of Corollary \ref{order}}]
         Consider any $(q+1)$-uniform family $\cT\subseteq \binom{\mathbb{N}}{q+1}$ with covering number $\tau(\cT)=q+1$. Let $\cF=\cF_{\cT,n}=\{T\cup G:T\in \cT, T\cap G=\emptyset, G\in \binom{[n]}{k-q-1}\}$. We claim that $\beta_{p,q}(\cF)\ge \binom{n-p-2q-1}{k-p-q-1}$ holds. Let $P,Q$ of sizes $p$ and $q$ be disjoint subsets of $[n]$. As $\tau(\cT)=q+1$, there exists $T\in \cT$ with $T\cap Q=\emptyset$. Since $k>p+q$, we have $|T\cup P|\le k$ and $T\cup P$ is disjoint with $Q$. By definition of $\cF$, all $k$-sets containing $T\cup P$ that are disjoint with $Q$ belong to $\cF$ and if $n\ge k+q$, then the number of such sets is $\binom{n-p-2q-1}{k-p-q-1}$. 

        On the other hand, by Theorem \ref{tau}, if $\tau(\cF)\le q$, then $\beta_{p,q}(\cF)=0$, otherwise $\beta_{p,q}(\cF)=O(n^{k-p-q-1})$.
\end{proof}

\begin{proof}[\textbf{Proof of Theorem \ref{0q}}]
 Let $\cF_0$ be a $(q+1)$-uniform intersecting family that achieves $\beta_{0,q}(\cF)=\beta(q)$. For $p,k,n$ with $k>p+k$, set $\cF=\cF_{p,k,n}=\cup_{F\in \cF_0}\nabla^n_k(F)$. Consider any $q$-subset $Q$ and $p$-subset $P$ of $[n]$ with $|P\cap Q|=0$. By definition of $\beta(q)$, there exist at least $\beta(q)$ sets $F_1,F_2,\dots,F_{\beta(q)}$ in $\cF_0$ that are disjoint with $Q$. If for some $i$, the set $F_i$ intersects $P$, then the number of sets containing $F_i\cup P$ and disjoint with $Q$ is $\binom{n-p-2q-1+|F_i\cap P|}{k-p-q-1+|P\cap F_i|}=\omega(n^{k-p-q-1+|F_i\cap P|})=\omega(n^{k-p-q})$ and all these sets belong to $\cF(P,\overline{Q})$. If that is not the case, then 
 \begin{equation*}
 \begin{split}
 |\cF(P,\overline{Q})| & \ge   \sum_{i=1}^{\beta(q)}|\nabla^n_k(F_i\cup P)(\emptyset,Q)|-\sum_{1\le i<j\le \beta(q)}|\nabla^n_k(F_i\cup F_j\cup P)(\emptyset,Q)| \\
 & \ge \beta(q)\binom{n-2q-p-1}{k-q-p-1}-O(n^{k-q-p-2}).     \end{split}
 \end{equation*}
This shows $\beta_{p,q}(k,n)\ge \beta(q)\binom{n-2q-p-1}{k-q-p-1}-O(n^{k-q-p-2})$.

Let $\cF\subseteq \binom{[n]}{k}$ be an intersecting family. If $\tau(\cF)\le q$, then by Theorem \ref{tau}, we have $\beta_{p,q}(\cF)=0$. So suppose $\tau(\cF)\ge q+1$, and let $\cT=\{T\in \cT_\cF: |T|=q+1\}, \cT^+=\cT_\cF\setminus \cT$. Note that by Theorem \ref{cover}, $|\cup_{T\in\cT_\cF}T|\le k^{k+1}$ so if $n$ is large enough we can pick $P$ to be disjoint with all $T\in \cT_\cF$. Then for any $T\in \cT^+$, we have $|\nabla^n_k(T\cup P)|\le \binom{n-p-q-2}{k-p-q-2}$. As $\cT_\cF$ is intersecting, so is $\cT$. Therefore we can pick $Q$ such that $\beta(q)\ge \cT(\emptyset,\overline{Q})$, i.e. the number of sets in $\cT$ that are disjoint with $Q$ is at most $\beta(q)$. Now, $\cF(P,\overline{Q})\subset \cup_{T\in \cT_\cF:Q\cap T=\emptyset}\nabla^n_k(T\cup P)$ and so $\beta_{p,q}(\cF)\le |\cF(P,\overline{Q})|\le \beta(q)\binom{n-p-q-1}{k-p-q-1}+(\binom{n-p-q-2}{k-p-q-2})$. This completes the proof of the theorem.
\end{proof}

For the proof of Theorem \ref{02}, we will need the following result on $\beta(2)$.

\begin{lemma}\label{charact}
    $\beta_{0,2}(n,3)=2$. Furthermore, if $\cF$ is a 3-uniform intersecting family with $\beta_{0,2}(\cF)=2$, then either $\cF$ is the Fano plane, or $\cF$ consists of 10 triples, one from every complement pair of 3-subsets of a six-element set. In both cases $\cT_{\cF}=\cF$.
\end{lemma}

\begin{proof}
    Let $\cF$ be a 3-uniform intersecting family with $\beta_{0,2}(\cF)=\beta(2)\ge 2$ (this last inequality follows from the Fano plane of which the $(0,2)$-dömdödöm is 2). By Theorem \ref{tau}, we know that $\tau(\cF)=3$. Hanson and Toft \cite{HT} proved that every 3-uniform intersecting family with $\tau(\cF)=3$ lives on at most 7 vertices. If the number of vertices in $\cF$ is at most 4, then $\tau(\cF)\le 2$. Next observe that $\beta_{0,2}(\binom{[5]}{3})=1$. 
    
    Suppose next $\cup_{F\in \cF}F=[6]$. Then $\beta_{0,2}(\cF)\ge 2$ implies that for every $x,y\in [6]$, $\cF$ contains at least 2 triples in $[6]\setminus \{x,y\}$. As every triple is contained in three 4-subsets of $[6]$, this implies $|\cF|\ge 10$. As $\cF$ is intersecting, we obtain that $|\cF|=10$, $\beta_{0,2}(\cF)=2$, and $\cF$ contains one triple from every complement pair of $[6]$. If $T\subseteq [6]$ with $|T|\ge 4$ was a minimal cover of $\cF$, then by minimality of $T$, $\cF$ would contain no triples in $T$, so no triples not containing elements of $[6]\setminus T$ contradicting $\beta_{0,2}(\cF)= 2$.

    Finally, assume $\cup_{F\in \cF}F=[7]$. If $\cF$ is maximal and $|F\cap F'|=1$ for all $F,F'\in \cF$, then any pair $x,y\in [7]$ is contained in at most one triple of $\cF$, so $|\cF|\le 7$. Dow, Drake, Füredi, and Larson showed \cite{DDFL} that every maximal 3-uniform intersecting family with covering number 3 contains at least 7 triples, so $|\cF|=7$ and every pair $x,y$ is contained in exactly one triple and so $\cF$ satisfies the axioms of projective planes, so $\cF$ is the Fano plane which possesses the property $\cT_\cF=\cF$. So assume $\cF$ contains $F,F'$ with $|F\cap F'|=2$, say $F=\{1,2,3\},F'=\{1,2,4\}$. Then $\beta_{0,2}(\cF)\ge 2$ implies the existence of $F_1,F_2\in \cF$ with $1,2\notin F_1,F_2$. Then by the intersecting  property of $\cF$, $3,4\in F_1\cap F_2$, so we can assume $F_1=\{3,4,5\},F_2=\{3,4,6\}$. Now $\beta_{0,2}(\cF)\ge 2$ implies the existence of $F_3,F_4\in \cF$ with $3,4\notin F_3,F_4$. By the intersecting property of $\cF$, they must be $\{5,6,1\}$ and $\{5,6,2\}$. But the covering number of these six triples is already 3, so $\cup_{F\in \cF}F=[6]$. This contradiction shows that we listed all possibilities.
\end{proof}

\begin{proof}[\textbf{Proof of Theorem \ref{02}}]
    The lower bound of (1) is given by the family $\cF^{2,3}_{n,k}=\{F\in \binom{[n]}{k}:|F\cap [3]|\ge 2\}$. Clearly, $\beta_{p,1}(\cF^{2,3}_{n,k})$ is achieved by $\cF^{2,3}_{n,k}(P,\overline{j})$ for any $j\in [3]$ and $p$-set $P$ disjoint with $[3]$ and $\beta_{p,1}(\cF^{2,3}_{n,k})=\cF^{2,3}_{n,k}(P,\overline{j})=\binom{n-3-p}{k-2-p}$.

    The lower bound of (2) is given by the Fano plane. More precesily, let $$\cF^{ano}=\{\{1,2,3\},\{3,4,5\},\{5,6,1\},\{2,4,6\},\{1,4,7\},\{3,6,7\},\{2,5,7\}\}$$ and for any $n,k$ let $\cF^{ano}_{n,k}=\cup_{F\in \cF^{ano}}\nabla^n_k(F)$. It is easy to see that $\beta_{p,2}(\cF^{ano}_{n,k})=|\cF^{ano}_{n,k}(P,\overline{Q})|$ for $p$-sets $P$ dosjoint with $[7]$ and pairs $Q\subseteq [7]$. For any such $Q$ there are exactly two $F,F'\in \cF^{ano}$ disjoint with $Q$, and their union $F\cup F'$ has size 5. Therefore $|\cF^{ano}_{n,k}(P,\overline{Q})|=2\binom{n-5-p}{k-3-p}-\binom{n-7-p}{k-5-p}$.

    \smallskip

    The proofs of the upper bounds are along the lines of that of Theorem \ref{order}. First we prove the upper bound of (1). Let $\cF\subseteq \binom{[n]}{k}$ be an intersecting family. Theorem \ref{tau} implies that if $\tau(\cF)=1$, then $\beta_{p,1}(\cF)=0$, while if $\tau(\cF)\ge 3$, then $\beta_{p,1}(\cF)=O(n^{k-3-p})=o(\binom{n-3-p}{k-2-p})$. So we can assume $\tau(\cF)=2$. We partition $\cT_\cF$ into $\cT_0\cup \cT^+$ with $\cT_0=\{T\in \cT_\cF:|T|=2\}$ and $\cT^+=\cT_\cF\setminus\cT_0$. If $\tau(\cT_0)=1$, then consider $j$ a shared element of all sets of $\cT_0$ and a $p$-set $P$ that is disjoint with all $T\in \cT_\cF$. Then $\cF(P,\overline{j})\subseteq \cup_{T\in \cT^+}\nabla^n_k(T\cup P)$ so, by Theorem \ref{cover}, its size and $\beta_{p,1}(\cF)$ is $O(n^{k-3-p})=o(\binom{n-3-p}{k-2-p})$. So we can assume $\tau(\cT_0)=2$ (as $\cT_0$ is intersecting). There is only one such family: $\cT_0=\{\{1,2\},\{2,3\},\{1,3\}\}$. Also, $\cT^+=\emptyset$,  since $\cT$ is intersecting and the only sets meeting all sets of $\cT_0$ are the ones in $\cT_0$. So for some $p$-set $P\subseteq [n]\setminus [3]$, we have $\beta_{p,1}(\cF)\le |\cF(P,\overline{1})|=|\{F\in \cF: \{2,3\}\cup P\subset F, 1\notin F\}|\le \binom{n-3-p}{k-2-p}$ as claimed.

    \smallskip
    
    To see the upper bound of (2), let $\cF\subseteq \binom{[n]}{k}$ be an intersecting family. Theorem \ref{tau} implies that if $\tau(\cF)\le 2$, then $\beta_{p,2}(\cF)=0$, while if $\tau(\cF)\ge 4$, then $\beta_{p,2}(\cF)=O(n^{k-4-p})=o(2\binom{n-5-p}{k-3-p}-\binom{n-7-p}{k-5-p})$. So we can assume $\tau(\cF)=3$. We again partition $\cT_\cF$ into $\cT_0\cup \cT^+$ with $\cT_0=\{T\in \cT_\cF:|T|=3\}$ and $\cT^+=\cT_\cF\setminus\cT_0$. If $\tau(\cT_0)\le 2$, then consider a 2-set $Q$ that is a cover of $\cT_0$ and a $p$-set $P$ that is disjoint with all $T\in \cT_\cF$. Then $\cF(P,\overline{Q})\subseteq \cup_{T\in \cT^+}\nabla^n_k(T\cup P)$ so its size and $\beta_{p,2}(\cF)$ is $O(n^{k-4-p})=o(2\binom{n-5-p}{k-3-p}-\binom{n-7-p}{k-5-p})$. So we can assume $\tau(\cT_0)=3$ (as $\cT_0$ is intersecting). If $\beta_{0,2}(\cT_0)=1$, then let $Q$ be a pair with $\cT_0(\emptyset,\overline{Q})=\{T_0\}$, and $P$ again a $p$-set disjoint with all $T\in \cT_\cF$. Then $\cF(P,\overline{Q})\subseteq \nabla^n_k(T_0\cup P)\cup \bigcup_{T\in \cT^+}\nabla^n_k(T\cup P)$ and so its size is at most $\binom{n-5-p}{k-3-p}+O(n^{k-4-p})<2\binom{n-5-p}{k-3-p}-\binom{n-7-p}{k-5-p}$ if $n$ is large enough. 

    Finally, if $\beta_{0,2}(\cT_0)=2$, then Lemma \ref{charact} implies that $\cT_0$ is either the Fano plane or contains one of every compelement pairs of triples of a six-element set $X$. Also in both cases, $\cT_{\cT_0}=\cT_0$. This implies $\cT^+=\emptyset$ as $\cT_\cF$ is intersecting, so every $T\in \cT^+$ contains a minimal cover $T'$ of $\cT_0$ which must be a member of $\cT_0$, so $T'$ cannot be a \textit{minimal} cover of $\cF$. Fix a 2-set $Q$ such that $\cT_0(\emptyset,\overline{Q})=\{T_1,T_2\}$ and a $p$-set $P$ disjoint with the vertices of the Fano plane or $X$. In the former case, we have $|T_1\cup T_2|=5$, while in the latter case $|T_1\cup T_2|=4$. Therefore, we have $|\cF(P,\overline{Q})|=2\binom{n-5-p}{k-3-p}-\binom{n-7-p}{k-5-p}$ or $|\cF(P,\overline{Q})|=2\binom{n-5-p}{k-3-p}-\binom{n-6-p}{k-6-p}<2\binom{n-5-p}{k-3-p}-\binom{n-7-p}{k-5-p}$. From the proof it follows that among \textit{maximal} intersecting families only $\cF^{ano}_{n,k}$ attains $\beta_{p,2}(\cF^{ano}_{n,k})=2\binom{n-5-p}{k-3-p}-\binom{n-7-p}{k-5-p}$. It is easy to see that a subfamily $\cF$  of $\cF^{ano}_{n,k}$ satisfies $\beta_{p,2}(\cF)=2\binom{n-5-p}{k-3-p}-\binom{n-7-p}{k-5-p}$ if and only if $\cF$ contains all sets in $\cF^{ano}_{n,k}$ that are supersets of exactly one set in $\cF^{ano}$.
\end{proof}

\section{Remarks}

In some proofs throughout the literature, instead of $\cF(P,\overline{Q})$, the family $\cF'(P,\overline{Q})=\{F\in \cF: P\cap F\neq \emptyset, Q\cap F=\emptyset\}$ is used, and thus one can introduce $$\beta'_{p,q}(\cF)=\min\{|\cF'(P,\overline{Q})|:|P|=p,|Q|=q, P\cap Q=\emptyset, P,Q\subseteq [n]\}$$ and $\beta'_{p,q}(n,k)=\max\{\beta'(\cF): \cF\subseteq \binom{[n]}{k}, \cF ~\text{is intersecting}\}$. By definition, it is clear that $\beta_{1,q}(n,k)\le \beta'_{p,q}(n,k)$ holds. Moreover, a similar proof to that of Theorem \ref{0q} shows that $\beta'_{p,q}(n,k)=(p+o(1))\beta_{1,q}(n,k)$.

Based on Theorem \ref{0q}, the main open problem is to determine $\beta(q)$. The maximum size of a $(q+1)$-uniform intersecting family $\cF$ with $\tau(\cF)=q+1$ is a much studied yet undetermined function which is clearly an upper bound on $\beta(q)$. Another simple inequality for $(q+1)$-uniform intersecting families $\cF$ with $\tau(\cF)=q+1$ is $\beta_{0,q}(\cF)\le \delta(\cF)-1$. Here for any vertex $x$ of minimum degree, one takes $Q=F\setminus \{x\}$ for any $x\in F\in \cF$.

\end{document}